\documentstyle[11pt]{article}

\title{NONSTANDARD DIGRAPHS}
\author{A.H. Zemanian}
\date{}

\begin{document}
\newcommand{\be}{\begin{equation}}
\newcommand{\ee}{\end{equation}}
\newcommand{\la}{\leftarrow}
\newcommand{\ra}{\rightarrow}
\newcommand{\hla}{\hookleftarrow}
\newcommand{\hra}{\hookrightarrow}
\newcommand{\dv}{\dashv}
\newcommand{\vd}{\vdash}
\newcommand{\as}{\asymp}
\newcommand{\lla}{\langle}
\newcommand{\rra}{\rangle}
\newcommand{\N}{I \kern -4.5pt N}
\newcommand{\sss}{^{*}\!}

\maketitle \baselineskip21pt

{\ Abstract --- Nonstandard graphs have been defined and examined in
prior works.  The present work does the same for nonstandard
digraphs.  Since digraphs have more structure than do graphs, the
present discussion requires more complicated definitions and yields
a variety of results peculiar to nonstandard digraphs.  A
nonstandard digraph can be obtained by means of an ultrapower
construction based on a sequence of digraphs or more elegantly by
using the transfer principle.  We use either or both techniques in
particular circumstances.  As special cases, we have the enlargement
of a single infinite digraph and also hyperfinite digraphs based on
sequences of finite digraphs.  Also examined are such ideas as
incidences and adjacencies for nonstandard arcs and vertices,
connectedness, components, and galaxies in nonstandard digraphs.

}

\section{Introduction}

Transfinite and nonstandard\footnote{It is preferable here to use
the adjective "nonstandard" instead of the prefix "hyper", commonly
used in nonstandard analysis,  because the word "hypergraph"
designates an entirely different kind of graph \cite{be}.  Also, we
will be using "standard" as a special case of "nonstandard".} graphs
have been constructed and examined in prior works. More recent works
on these subjects are \cite{B7}, \cite{B8}, and \cite{gal}.
However, transfinite and nonstandard generalizations of digraphs had
not yet been achieved.   We now aim to do so.  In this present
report, we discuss nonstandard digraphs.  In two subsequent reports,
we will investigate transfinite digraphs and digraphs that are both
transfinite and nonstandard.  Since digraphs have more structure
than do graphs, our present discussions, though similar to those for
graphs, require more complicated definitions and yield more detailed
and broader results.

Our notations and symbols are the same as those specified in
\cite[Section 1.1]{B8}.  Let us mention here some of them.  Braces
$\{\ldots\}$ denote a set;  its elements within the braces are all
distinct, and the order in which they are listed is not significant.
On the other hand, angle brackets $\lla\ldots\rra$ denote a
sequence;  its elements within the angle brackets have an imposed
order from left to right;  also, those elements may repeat.  As a
special case, we have an ordered pair, that is, a two-element
sequence $\lla a,b\rra$ with $a$ preceding $b$.  $\N$ denotes the
set of natural numbers: $\{0,1, 2, \ldots \}$.  $\cal F$ denotes a
nonprincipal ultrafilter on  $\N$, which will remain fixed
throughout this work.

\section{Standard digraphs}

Before presenting the definitions for nonstandard digraphs, let us
state explicitly what we mean by a "standard digraph."  We will use
a rather different, but virtually equivalent, definition of a
standard digraph as compared to the conventional
definition.\footnote{See, for instance, \cite[Section 22]{wi}.}  The
reason for this is that we wish to construct standard digraphs in
the same way as will be done for transfinite digraphs.

We start with  a set $A$ of arcs, where each arc $a\in A$ is an
ordered pair $a=\lla s,t\rra$ of {\em ditips}.  We refer to $s$ as
the {\em intip} of $a$ and to $t$ as the {\em outtip} of $a$, and we
view $a$ as having a direction from $s$ to $t$.  Then, the union of
all the arcs is a set $T$ of ditips with each intip having a
corresponding outtip in $T$ in accordance with the arcs.  We will
also let $T_{i}$ (resp. $T_{o}$) denote the set of intips (resp.
outtips), and thus $T=T_{i}\cup T_{o}$.

We now partition $T$ arbitrarily.  Each set $v$ in that partition is
a {\em vertex}.  $V$ will denote the set of vertices.  Accordingly,
we say that each arc $a=\lla s,t\rra$ is {\em directed} from the
vertex $u$ containing its intip $s$ toward the vertex $v$ containing
its outtip $t$.  Possibly, $u$ and $v$ are the same vertex,
resulting in the arc being a {\em self-loop}.  We say that $u$ and
$a$ are {\em incident inward} and that $v$ and $a$ are {\em incident
outward}.\footnote {In order to relate our unusual version of a
digraph to the conventional one, we can at this point identify an
arc $a$ as the ordered pair $a=\lla u,v\rra$.  However, our present
construction allows self-loops and parallel arcs---in contrast to
the conventional case.}

Finally, a (standard) {\em digraph} $D$ is the pair \be D=\{A,V\}.
\label{2.1} \ee

The "underlying graph" $G$ of $D$ is obtained by removing the
directions of the arcs.  Thus, each intip and each outtip becomes
merely a {\em tip} with no implication of a direction.  An arc
$a=\lla s,t\rra$ becomes a {\em branch} $b=\{ s,t\}$ with $s,t\in
T$, where $T$ is the set of tips.  The set $A$ of arcs becomes a set
$B$ of branches.  Moreover, the partition of the set $T$ of tips now
becomes a set $X$ of {\em nodes}.  Finally, the {\em underlying
graph} $G$ of $D=\{A,V\}$ is \be G=\{B,X\}            \label{2.2}
\ee This structure is discussed in more detail in \cite[Section
2.2]{B8}.

In subsequent parts of this three-part work, we will be dealing with
transfinite digraphs.  These appear in a hierarchy of
transfiniteness ranked by the natural numbers and subsequently by
the countable ordinals.  In that case, the standard digraphs and its
vertices will be assigned the rank 0, and the notation (\ref{2.1})
will be replaced by
\[ D^{0}\;=\;\{ A,V^{0}\}. \]
Also, the intips and outtips will be called $(-1)$-intips and
$(-1)$-outtips, respectively.  Furthermore,    $T^{-1}$ will denote
the set of  the ditips of all the arcs, with the notations $s$ and
$t$ replaced by $s^{-1}$ and $t^{-1}$, respectively.

\section{Nonstandard digraphs}

Let $\lla D_{n}:n\in \N\rra$ be some chosen and fixed sequence of
standard digraphs.  Here, $D_{n}=\{A_{n},V_{n}\}$ with $A_{n}$ and
$V_{n}$ being respectively the set of arcs and the set of vertices
for $D_{n}$.  These digraphs may overlap;  that is, for $n\not = m $
we may have            $A_{n}\cap A_{m}\not=\emptyset$.\footnote{In
fact, we can view the $D_{n}$ as being subgraphs of a large digraph
$D=\{A,V\}$, where $A$ contains $\cup_{n=0}^{\infty}A_{n}$ .  But,
this idea will not be pursued.}  In addition, $\cal F$ will denote a
nonprincipal ultrafilter.  It will be understood henceforth that
$\cal F$ is chosen and fixed.

Next, let $\lla a_{n}\rra = \lla a_{n}: n\in\N\rra$ be a sequence of
arcs with $a_{n}\in A_{n}$ for every $n\in\N$. A {\em nonstandard
arc} $\bf a$ is an equivalence class of all such sequences of arcs.
By "equivalence" we mean that every two such sequences $\lla
a_{n}\rra$ and $\lla \acute{a}_{n}\rra$ are taken to be {\em
equivalent} if $\{n:a_{n}=\acute{a}_{n}\}\in{\cal F}$.  In this
case, we write "$\lla a_{n}\rra=\lla \acute{a}_{n}\rra$ a.e." or say
that $a_{n}=\acute{a}_{n}$ for almost all $n$.  We also use the
notation ${\bf a}=[a_{n}]$ where the $a_{n}$ are the members of one
i.e., any one) of the sequences in the equivalence class.

That this truly partitions the set of all such sequences is seen as
follows.  Reflexivity and symmetry are obvious.  As for
transitivity, let $\lla a_{n}\rra$, $\lla \acute{a}_{n}\rra$, and
$\lla \tilde{a}_{n}\rra$ be three such sequences with $\lla
a_{n}\rra$ and $\lla \acute{a}_{n}\rra$ being equivalent and $\lla
\acute{a}_{n}\rra$ and $\lla\tilde{a}_{n}\rra$ being equivalent .
Then, $N_{a\acute{a}}=\{n:a_{n}=\acute{a}_{n}\}\in{\cal F}$,
$N_{\acute{a}\tilde{a}}=\{n:\acute{a}_{n}=a\tilde{a}\}\in{\cal F}$,
and $N_{a\tilde{a}}=\{n:a_{n}=\tilde{a}\}\supseteq
N_{a\acute{a}}\cap N_{\acute{a}\tilde{a}}\in{\cal F}$.  Hence,
$N_{a\tilde{a}}\in{\cal F}$, which asserts that $\lla a_{n}\rra$ and
$\lla \tilde{a}_{n}\rra$ are equivalent.

Whenever $\lla a_{n}\rra$ and $\lla \acute{a}_{n}\rra$ are
equivalent, their corresponding sequences $\lla s_{n}\rra$ and $\lla
\acute{s}_{n}$ of intips are perforce equivalent, too, because each
intip uniquely determines its arc, and conversely.  Thus,
$\{n:s_{n}=\acute{s}_{n}\}=\{n:a_{n}=\acute{a}_{n}\}\in{\cal F}$.
For the same reason, the corresponding sequences of outtips ,
namely, $\lla t_{n}\rra$ and $\lla \acute{t}_{n}\rra$ are
equivalent, too.

So far, we have, from the above partition of the set of sequences
$\lla a_{n}\rra$ of arcs, a set  $^{*}\!A$ of nonstandard arcs ${\bf
a}=[a_{n}]$.  Correspondingly, we get a partition of the set of
sequences $\lla s_{n}\rra$ of intips (resp. a partition of the set
of sequences $\lla t_{n}\rra$ of outtips), and the sets of that
partition are the {\em nonstandard intips} (resp. {\em nonstandard
outtips}).  We let ${\bf s}=[s_{n}]$ (resp. ${\bf t}=[t_{n}]$) be a
typical nonstandard intip (resp. nonstandard outtip), and we then
have the nonstandard arc ${\bf a}=\lla{\bf s}, {\bf t}\rra$.  We let
$T_{i}$ (resp. $T_{o}$) denote the set of all nonstandard intips
(resp. the set of all nonstandard outtips).  We now wish to
construct the set $^{*}\! V$ of all "nonstandard vertices."

For each $n\in\N$, let $p_{n}$ be a tip of an arc in $A_{n}$;  that
is, $p_{n}$ is either an intip or an outtip of that arc.  Then,
consider the sequence $\lla p_{n}\rra$.  Let $N_{i}$ be the set of
all $n$ for which $p_{n}$ is an intip, and let $N_{o}$ be the set of
all $n$ for which $p_{n}$ is an outtip.  Thus, $N_{i}\cap
\N_{o}=\emptyset$ and $N_{i}\cup N_{o}=\N$.  So, either $N_{i}$ or
$N_{o}$ (but not both) is a member of  $\cal F$.  If it is $N_{i}$
(resp. $N_{o}$), we can show that $\lla p_{n}\rra$ is the
representative of a nonstandard intip (resp. a nonstandard outtip)
as follows.

Let $\lla p_{n}\rra$ and $\lla q_{n}\rra$ be two equivalent
sequences of ditips.  Remember that these are taken to be equivalent
if $p_{n}=q_{n}$ for almost all $n$. This equivalence partitions the
set of all sequences of ditips into equivalence classes. Indeed,
reflexivity and symmetry are obvious, and transitivity follows as
usual.  Each such equivalence class is taken to be a {\em
nonstandard ditip} ${\bf p}=[p_{n}]$.  Moreover, if $[p_{n}]$ is a
nonstandard intip, then so, too, must be $[q_{n}]$.  Indeed, we have
$N_{pq}=\{n:p_{n}=q_{n}\}\in{\cal F}$.  Moreover, the set $N_{p}$ of
all $n$ for which $p_{n}$ is an intip is also a member of $\cal F$.
Let $N_{q}$ be the set of all $n$ for which $q_{n}$ is an intip.
Now, $N_{q}\supseteq N_{p}\cap N_{pq}\in{\cal F}$.  Hence,
$N_{q}\in{\cal F}$.  Thus, $[q_{n}]$ is a nonstandard intip.  In the
same way, it follows that, if  $[p_{n}]$ is a nonstandard outtip,
then so, too, is $[q_{n}]$.

More notation:  Let $p_{n}$ and $q_{n}$ be two ditips of $D_{n}$,
not necessarily of the same kind. That is, one may be an intip and
the other an outtip or they may be both intips or both outtips.  If
$p_{n}$ and $q_{n}$ are members of the same vertex in $V_{n}$ (resp.
in different vertices of $V_{n}$), we say that $p_{n}$ and $q_{n}$
are {\em shorted together} (resp. {\em not shorted together}), and
we write $p_{n}\as q_{n}$ (resp. $p_{n}\not\as q_{n}$).

Next, let ${\bf p}=[p_{n}]$ and ${\bf q}=[q_{n}]$ be two nonstandard
ditips.  Let $N_{pq}=\{n:p_{n}\as q_{n}\}$ and $N_{pq}^{c}=\{n:
p_{n}\not\as q_{n}\}$.  Either $N_{pq}\in{\cal F}$ or
$N_{pq}^{c}\in{\cal F}$, but not both.  If $N_{pq}\in{\cal F}$
(resp. $N_{pq}^{c}\in{\cal F}$), we say that $\bf p$ and $\bf q$ are
{\em shorted together} and we write ${\bf p}\as{\bf q}$ (resp. $\bf
p$ and $\bf q$ are {\em not shorted together} and we write ${\bf
p}\not\as{\bf q}$).  Furthermore, we take it that $\bf p$ is shorted
to itself: ${\bf p}\as{\bf p}$.  This shorting is an equivalence
relation for the set of all nonstandard ditips.  Indeed, with
reflexivity and symmetry again being obvious, consider transitivity.
Assume ${\bf p}\as{\bf q}$ and ${\bf q}\as{\bf r}$.  Since $\{n:
p_{n}\as q_{n}\}\cap\{n:q_{n}\as r_{n}\}\subseteq\{n: p_{n}\as
r_{n}\}$, it follows that ${\bf p}\as{\bf r}$.  The resulting
equivalence classes are the {\em nonstandard vertices}.

This definition is independent if the representative sequences
chosen for the nonstandard ditips.  Indeed, let ${\bf
p}=[p_{n}]=[\tilde{p}_{n}]$ and ${\bf q}=[q_{n}]=[\tilde{q}_{n}]$.
Set $N_{p}=\{n:p_{n}=\tilde{p}_{n}\}\in{\cal F}$ and
$N_{q}=\{n:q_{n}=\tilde{q}_{n}\}\in{\cal F}$. Assume
$[p_{n}]\as[q_{n}]$.  Thus, $N_{pq}=\{n:p_{n}\as q_{n}\}\in{\cal
F}$.  We want to show that $N_{\tilde{p}\tilde{q}}\}$ is a member of
$\cal F$, so that $[\tilde{p}_{n}]\as[\tilde{q}_{n}]$.  We have
$(N_{p}\cap N_{q}\cap N_{pq})\subseteq N_{\tilde{p}\tilde{q}}$.
Hence, $N_{\tilde{p}\tilde{q}}=\{n:\tilde{p}\as\tilde{q}\}\in {\cal
F}$, whence our conclusion.

Altogether, we have defined a nonstandard vertex $\bf v$ to some set
in the partition of the set of nonstandard ditips induced by the
shorting $\as$.  $^{*}\! V$ will denote the set of nonstandard
vertices.

Thus, we now have the {\em nonstandard digraph}\footnote{To conform
with common terminology in nonstandard analysis, we could have
called this a "hypergraph" but will not do so because that term is
used for and entirely different kind of graph \cite{be}.} \be ^{*}\!
D\;=\;\{^{*}\! A, ^{*}\! V\}  \label{3.1} \ee obtained from the
given sequence $\lla D_{n}\rra$ of standard digraphs.

Let us note that we obtained $^{*}\! D$ by starting from a given
sequence $\lla D_{n}:n\in \N\rra$, where $D_{n}=\{A_{n},V_{n}\}$.
But, any other sequence $\lla \tilde{D}_{n}:n\in\N\rra$, where
$\tilde{D}_{n}=\{\tilde{A}_{n},\tilde{V}_{n}\}$ could be used to get
the same $^{*}\! D$ so long as $\{n:A_{n}=\tilde{A}_{n}\}$ and
$\{n:V_{n}=\tilde{V}_{n}\}$ are both members of $\cal F$.  In this
regard, see \cite[Theorem 12.1.1]{go}.

Finally, the {\em underlying nonstandard graph} $^{*}\! G$ of
$^{*}\! D$ is obtained simply by removing the directions of the arcs
in each $D_{n}$ to get branches in a standard graph $G_{n}$.  When
doing this, oppositely directed arcs incident to the same two
vertices become parallel branches, but we allow parallel branches.
Then, the equivalence partitioning of the set of sequences of arcs
becomes an equivalence partitioning of the set of sequences of
branches to yield the nonstandard branches, the set of which is
denoted by $^{*}\!B$. In the same way, the ditips of an arc in
$D_{n}$ become the tips of a branch in $G_{n}$.  Then, the
corresponding equivalence partitioning of the set of all sequences
of ditips that yielded the nonstandard ditips becomes an equivalence
partitioning of the set of all sequences of tips, and the sets of
that partitioning become the nonstandard tips. Next, the shortings
of the nonstandard ditips to get the nonstandard vertices is
mimicked by shortings of nonstandard tips to get the nonstandard
nodes, the set of which is denoted by $^{*}\!X$. Finally, the
underlying graph of $^{*}\!D=\{ ^{*}\! A,^{*}\! V\}$ is
\[ ^{*}\!G\;=\; \{^{*}\!B,^{*}\!X\}. \]

\section{Special cases}

One special case arises when all the $D_{n}$ are the same standard
graph $D=\{A,V\}$.  In this case, $^{*}\!D=\{^{*}\!A,^{*}\!V\}$ will
be called an {\em enlargement} of $D$.

If $D$ is a finite digraph, each arc ${\bf a}\in\, ^{*}\!A$ can be
identified with an arc $a\in\, A$ because the enlargement of a
finite set is the set itself.  In this case, every nonstandard
vertex ${\bf v}\in\, ^{*}\! V$ can be identified with a vertex $v\in
V$. Thus, $^{*}\! D=D$.

On the other hand, if $D=\{A,V\}$ is a conventionally infinite
digraph, that is, if $A$ is an infinite set, then the enlargement
$\sss A$ of A has more elements than $A$, namely, nonstandard arcs
that are different from the standard arcs (i.e., $\sss A\backslash
A$ is not empty).  However, $V$ may or may not be an infinite set.
If $V$ is infinite, $\sss V\backslash V$ is not empty, too.  But, if
$V$ is finite, then each vertex in $\sss V$ can be identified with a
vertex in $V$. In either case, $\sss D$ is a proper enlargement of
$D$ since $\sss A$ is a proper enlargement of $A$.  This special
case is examined again in Section 10.

Another special case arises when almost all of the
$D_{n}=\{A_{n},V_{n}\}$ are (possibly different)  finite digraphs.
We now call the nonstandard digraph $\sss D=\{\sss A,\sss V\}$ a
{\em hyperfinite digraph}.  As a result, we can lift many theorems
concerning finite digraphs into theorems about hyperfinite digraphs.
This can be done by writing the theorems about finite digraphs in
symbolic-logic notation and then applying the transfer principal.

\section{Incidences}

Given an arc $a=\lla s,t\rra$ and a vertex $u$ containing the intip
$s$ of $a$, we say that $a$ and $u$ are {\em incident  inward} and
write $a\hla u$ or $u\hra a$.  Similarly, if a vertex $v$ contains
the outtip $t$ of $a$, we say that $a$ and $v$ are {\em incident
outward} and write $a\hra v$ or $v\hla a$.\footnote{Our symbols
$\hra$ and $\hla$ are not the implication symbols $\rightarrow$ and
$\leftarrow$ used in symbolic languages.}  For the symbolic
sentences used below, we let $T_{i}$ (resp. $T_{o}$) be the set of
intips (resp. outtips). Then, $u\hra a$ will mean that $(\exists
s\in T_{i})(\exists u\in V)(\exists a\in A)(s\in u\wedge s\in a)$.
Also, $a\hra v$ will mean $(\exists t\in T_{o})(\exists v\in
V)(\exists a \in A)(t\in a\wedge t\in v)$.  Thus, "inward" and
"outward" express directions with respect to $a$ (not with respect
to $v$).  So, if the vertex $u$ contains the intip of $a$ and the
vertex $v$ contains the outtip of $a$, we may write $u\hra a\hra v$
or $v\hla a\hla u$. It is possible that $u$ and $v$ are the same
vertex, in which case $a$ is a {\em self-loop}.

We can express these incidences in symbolic notation as follows.  In
the standard case, we have \be (\exists a\in A)(\exists u\in
V)(\exists v\in V)(u\hra a\hra v)   \label{5.1} \ee Upon applying
the transfer principle, we obtain the symbolic sentence in terms of
nonstandard quantities: \be (\exists {\bf a}\in{\bf A})(\exists{\bf
u}\in{\bf V})(\exists {\bf v}\in{\bf V})({\bf  u}\hra{\bf a}\hra{\bf
v})  \label{5.2} \ee The words "incident inward"  and "incident
outward" then apply for ${\bf u}\hra{\bf a}$ and ${\bf a}\hra{\bf
v}$, respectively.

Alternatively, in terms of an ultrapower construction, we can take
it that, for almost all $n$, we have $a_{n}\in A_{n}$, $u_{n}\in
V_{n}$, $v_{n}\in V_{n}$, and then can require that $\{n: u_{n}\hra
a_{n}\hra v_{n}\}\in{\cal F}$ in order to obtain (\ref{5.2}) again.

\section{Adjacencies}

For a standard digraph $D=\{A,V\}$, two standard vertices $u,v\in V$
are called {\em adjacent} when the following is true:
\[ (\exists u,v\in V)(\exists a=\lla s,t\rra\in A)((s\in u\wedge t\in v)\vee(s\in v\wedge t\in v)).  \]
By transfer, we have {\em adjacency} for two nonstandard vertices
${\bf u},{\bf v}\in \sss V$ for $\sss D=\{\sss A,\sss V\}$ when the
following is true:
\[ (\exists {\bf u},{\bf v}\in\, \sss V)(\exists {\bf a}=\lla {\bf s},{\bf t}\rra\in\, \sss A)(({\bf s}\in{\bf u}\wedge{\bf t}\in{\bf v})\vee({\bf s}\in{\bf v}\wedge{\bf t}\in{\bf u})) \]

Similarly, two standard arcs $a,c\in A$ are called adjacent if the
following is true:
\[ (\exists w\in V)(\exists a=\lla s_{a},t_{a}\rra\in A)(\exists c=\lla s_{c},t_{c}\rra\in A)\\
((s_{a},s_{b}\in w)\vee( s_{a},t_{b}\in w)\vee(t_{a},s_{b}\in
w)\vee(t_{a},t_{b}\in w)) \] Again by transfer, we get the
definition for {\em adjacency} for two nonstandard arcs ${\bf
a},{\bf c}\in\, \sss A$ by using boldface notation for the vertices,
arcs, and tips and replacing $V$ and $A$ by $\sss V$ and $\sss A$,
respectively.

\section{Dipaths and semipaths}

A {\em standard finite dipath} $P$ in $D=\{A,V\}$ is defined as
follows:
\[ (\exists k\in\N\setminus\{0\})(\exists a_{0}, a_{1}, \ldots, a_{k-1}\in A)(\exists v_{0}, v_{1},\ldots, v_{k}\in
V) \]\be (v_{0}\hra a_{0}\hra v_{1}\hra a_{1}\hra v_{2}\hra
\cdots\hra v_{k-1}\hra a_{k-1}\hra v_{k})  \label{7.1} \ee It is
understood that all the arcs and vertices herein are distinct.  The
length $|P|$ of $P$ is the number of arcs herein; i.e., $|P|=k$.

By transfer, we get the following definition of a {\em nonstandard
hyperfinite dipath} $\sss P$ in $\sss D=\{\sss A,\sss V\}$:
\[ (\exists {\bf k}\in \N\setminus\{0\})(\exists {\bf a}_{0}, {\bf a}_{1}, \ldots ,{\bf a}_{{\bf k}-1}\in\sss A)(\exists {\bf v}_{0}, {\bf v}_{1}, \ldots ,{\bf v}_{\bf k}\in \sss V) \]
\be ({\bf v}_{0}\hra{\bf a}_{0}\hra{\bf v}_{1}\hra{\bf
a}_{1}\hra{\bf v}_{2} \ldots \hra{\bf v}_{{\bf k}-1}\hra {\bf
a}_{{\bf k}-1}\hra {\bf v}_{\bf k}) \label{7.2} \ee Because $\sss
D=[D_{n}]$, where the $D_{n}$ may be finite digraphs that grow
unlimitedly in size as $n$ increases through $\N$ or indeed may be
infinite digraphs, ${\bf k}$ may be a hypernatural number in $\sss
\N\setminus \N$.   The length $|\sss P|$ of $\sss P$ is $\bf k$.

The definition of a {\em nonstandard hyperfinite diloop} $\sss L$ is
obtained by transfer from the definition of a standard finite diloop
$L$;  that is, for $L$ we may use the definition (\ref{7.1}) with
the proviso that $v_{k}=v_{0}$.  Thus, by setting ${\bf v}_{\bf
k}={\bf v}_{0}$ in (\ref{7.2}), we obtain the definition of $\sss
L$.

A {\em finite semipath} $P_{s}$ in $D$ is defined as is a dipath
except that the directions of the arcs are ignored.  In particular,
the symbols $v\vd a$ and $a\dv$ will both mean that a tip in the
vertex $v$ is a member of the arc $a$.  That tip could be either an
intip or an outtip of $a$.  Similarly, the symbol $u \vd a\dv v$
means that a tip of $a$ is a member of $u$ and the other tip of $a$
is a member of $v$.  On the other hand, $a\dv v\vd b$ denotes that a
tip of $v$ is a member of $a$ and that another tip of $v$ is member
of the arc $b$.

Then, $P_{s}$ is defined by
\[ (\exists k\in\N\setminus\{0\})(\exists a_{0}, a_{1}, \ldots , a_{k-1}\in A)(\exists v_{0}, v_{1}, \ldots ,v_{k}\in V) \]
\be ( v_{0}\vd a_{0}\dv v_{1}\vd a_{1}\dv v_{2}\vd \ldots\ \dv
v_{k-1}\vd a_{k-1}\dv v_{k})  \label{7.3} \ee Here, too, it is
understood that all the arcs and vertices are distinct.  Actually,
$P_{s}$ can be identified as a path in the underlying graph $G$ of
$D$.  The length $|P_{s}|$ of $P_{s}$ is $k$.

By transfer, we get the definition of a {\em nonstandard hyperfinite
semipath} $\sss P_{s}$ in $\sss D=\{\sss A,\sss V\}$.  It can be
obtained from (\ref{7.3}) by writing boldface notation for $k$, $a$,
and $v$ and replacing $\N$, $A$, and $V$ by $\sss \N$, $\sss A$, and
$\sss V$. We have ${\bf k}\in \sss \N$ and possibly ${\bf k}\in\sss
\N\setminus\N$.  The length  of $\sss P_{s}$ is $|\sss P_{s}|={\bf
k}$.

The definition of a {\em standard finite diloop} is obtained from
(\ref{7.3}) by setting $v_{k}=v_{0}$.  For a {\em nonstandard
hyperfinite diloop}, make the same replacements as before.

\section{Connectedness and components}

The ideas of strong connectedness, unilateral connectedness, and
weak connectedness for standard digraphs are well-known \cite{hnc}
and need not be explicated here.  These ideas transfer directly to
the definition for connectedness between two vertices $\bf u$ and
$\bf v$ in the nonstandard digraph $\sss D=\{\sss A,\sss V\}$:

$\bf u$ and $\bf v$ are called {\em strongly connected} if there
exists a hyperfinite dipath from $\bf u$ to $\bf v$ and also a
hyperfinite dipath from $\bf v$ to $\bf u$.

$\bf u$ and $\bf v$ are called {\em unilaterally connected} if there
exists a hyperfinite dipath from one of those vertices to the other.

$\bf u$ and $\bf v$ are called {\em weakly connected} if there
exists a hyperfinite semipath terminating at $\bf u$ and $\bf v$.

Thus, strong connectedness implies unilateral connectedness, which
in turn implies weak connectedness.

Also, $\bf u$ and $\bf v$ are called {\em disconnected} if there is
no path of any kind between them.

A nonstandard digraph $\sss D=\{\sss A,\sss V\}$ is said to be {\em
strong} (resp. {\em unilateral}, resp. {\em weak}) if every two
vertices ${\bf u},{\bf v} \in\sss V$ are strongly connected (resp.
unilaterally connected, resp. weakly connected).  Also, $\sss D$ is
called {\em strictly unilateral} (resp. {\em strictly  weak}) if
$\sss D$ is unilateral but not strong (resp. $\sss D$ is weak but
not unilateral).  Moreover, $\sss D$ is called {\em disconnected} if
it has two vertices $\bf u$, $\bf v$ that are disconnected.

We turn now to the ideas of "subdigraphs" and "reduced digraphs".
Given the standard digraph $D=\{A,V\}$, let $A_{s}$ be a subset of
$A$, that is, $A_{s}$ consists of some but not necessarily all of
the arcs in $A$.  Furthermore, let $V_{s}$ be the subset of $V$
consisting of those vertices in $V$ having at least one tip (i.e.,
an intip or an outtip) belonging to an arc in $A_{s}$.  Then, the
{\em subdigraph $D_{s}$ of $D$ induced by $A_{s}$} is the doublet
\be D_{s}\;=\;\{A_{s},V_{s}\}           \label{8.1} \ee We can also
define a subdigraph $D_{ss}\;=\;\{A_{ss},V_{ss}\}$ of $D_{s}$ by
choosing a subset $A_{ss}$ of $A_{s}$ and defining $V_{ss}$ from
$A_{ss}$ as $V_{s}$ was defined from $A_{s}$.

Note that $D_{s}$ need not be a digraph by itself because $V_{s}$
may contain a vertex having a tip belonging to an arc not in
$A_{s}$.  To overcome this anomaly, we can {\em reduce} such
vertices as follows.  Reduce each such vertex $v$ in $V_{s}$ by
removing those tips in $v$ that are not members of arcs in $A_{s}$.
The remaining set of tips is the {\em reduced vertex} $v_{r}$;
$v_{r}$ is not empty.  The set of reduced vertices will be denoted
by $V_{r}$, and the pair $\{A_{s},V_{r}\}$ can be called the {\em
reduced digraph} $D_{r}$ of $D$ {\em induced by the chosen subset
$A_{s}$ of} $A$.  In the following, we will be dealing with the
subdigraph (\ref{8.1}) rather than the reduced digraph $D_{r}$.

The definition of a "nonstandard arc-induced subdigraph" $\sss
D_{s}=\{\sss A_{s},\sss V_{s}\}$ of $\sss D\;=\;\{\sss A,\sss V\}$
can be obtained by transfer of the above standard definition. Just
to vary our discussion, let us present the nonstandard definition by
means of an ultrapower construction of $\sss D$.  We start  with a
sequence $\lla D_{n}\rra=\lla A_{n},V_{n}\rra$ of standard digraphs,
where $n\in\N$, and also with a sequence $\lla
D_{s,n}\rra=\lla\{A_{s,n},V_{s,n}\}\rra$, where each $D_{s,n}$ is an
arc-induced subdigraph of $D_{n}$.  Let $[a_{n}]$ denote a
nonstandard arc in $\sss A$;  thus, $[a_{n}]$ is in $\sss A_{s}$ if
$\{n:a_{n}\in A_{s,n}\}\in{\cal F}$.  Furthermore, let $V_{s,n}$ be
the set of vertices in $V_{n}$ such that each vertex contains at
least one tip of an arc in $A_{s.n}$. Then, $[v_{n}]\in\,\sss V_{s}$
if $\{n:v_{n}\in V_{s,n}\}\in{\cal F}$. The intersection of these
two sets defining $\sss A_{s}$ and $\sss V_{s}$ is also in $\cal F$.
As a result, we obtain the subdigraph $\sss D_{s}=\{\sss A_{s},\sss
V_{s}\}$, which we call a {\em nonstandard arc-induced subdigraph
of} $\sss D$ (induced by the arcs in $\sss A_{s}$).  More concisely,
we refer to $\sss D_{s}$ as a {\em subdigraph of} $\sss D$.

We turn now to the concept of "components" in $\sss D=\{\sss A,\sss
V\}$. A {\em strong component} of $\sss D$ is a maximal set of
nonstandard vertices that are pairwise strongly connected.

By substituting "unilateral" or "weak" for "strong" in the preceding
paragraph, we get the definitions of {\em unilateral component} or
{\em weak component} in $\sss D$, respectively.  It follows that, in
$\sss D$, a strong component is a subset of a unilateral component,
and the latter is a subset of a weak component.

\section{Bounds on the number of arcs in a nonstandard hyperfinite digraph}

As was mentioned in Section 4, results concerning finite digraphs
can be extended directly to nonstandard  hyperfinite digraphs by
means of transfer.  As an example of this, let us transfer certain
bounds on the number $q$ of arcs of a standard digraph having $p$
vertices.

Let $\sss D_{f}$ be a nonstandard hyperfinite digraph with no
parallel arcs and no self-loops (i.e., for almost all of the finite
$D_{n}$ from which $\sss D_{f}$ is obtained, there are no parallel
arcs and no self-loops).  Also, let $\sss D_{f}$ have {\bf q} arcs
and {\bf p} vertices.  Here $\bf q$ and $\bf p$  are hypernatural
numbers.  We can obtain bounds on $\bf q$ in terms of $\bf p$ for
various categories of connectedness by transferring results on
standard digraphs, as stated for instance in \cite[pages
71-75]{hnc}. Specifically, we have the following:

If $\sss D_{f}$ is complete and symmetric (i.e., if $({\bf u},{\bf
v})$ and $({\bf v},{\bf u})$ are arcs in $\sss A$ for every pair of
vertices ${\bf u}, {\bf v}\in{\sss V})$, then ${\bf q}={\bf p}({\bf
p}-1)$.

If $\sss D_{f}$ is disconnected, then $0\leq{\bf q}\leq({\bf
p}-1)({\bf p}-2)$.

If $\sss D_{f}$ is strictly weak, then ${\bf p}-1\leq{\bf
q}\leq({\bf p}-1)({\bf p}-2)$ and ${\bf p}\geq 3$.

If $\sss D_{f}$ is strictly unilateral, then ${\bf p}-1\leq{\bf
q}\leq({\bf p}-1)^{2}$.

If $\sss D_{f}$ is strong and if ${\bf p}>1$, then ${\bf p}\leq{\bf
q}\leq{\bf p}({\bf p}-1)$.

\section{The galaxies of nonstandard enlargements of infinite digraphs}

The discussion in this section is much like that for enlargements of
graphs \cite{gal}, but there is more to say regarding enlargements
of digraphs.

We now start with a standard digraph $D=\{A,V\}$ having an infinity
of arcs and an infinity of vertices.  As always, we assume that $D$
is weakly connected.  Now, $\sss D=\{\sss A,\sss V\}$ denotes the
nonstandard enlargement of $D$.  Consequently, ${\sss A}\setminus A$
and ${\sss V}\setminus V$ are both nonempty and contain nonstandard
arcs and nonstandard vertices, respectively.

We define the "galaxies" of $\sss D$ in the same way as was done for
the enlargement of an infinite graph \cite[Section 3]{gal}.  Let us
be specific here.  The length $|P_{uv}|$ of any semipath $P_{uv}$
connecting two vertices $u$ and $v$ in $D$ is the number of arcs in
$P_{uv}$.  The distance $d(u,v)$ between $u$ and $v$ is $d(u,v)=
\min \{|P_{uv}|\}$, where the minimum is taken over all the
semipaths terminating at $u$ and $v$.  Also, for any vertex $u$, we
set $d(u,u)=0$.  $d$ is a metric for $V$, with reflexivity and
symmetry being obvious and the triangle inequality being readily
shown.   So, we can view $V$ as being a metric space with $d$ as its
metric.  Now, $d$ can be extended into an internal function $\bf d$
mapping the Cartesian product ${\sss V}\times{\sss V}$ into the set
of hypernaturals $\sss \N$, where, for any ${\bf u}=[u_{n}]$ and
${\bf v}=[v_{n}]$, ${\bf d}({\bf u},{\bf v})$ is defined by
\[ {\bf d}({\bf u},{\bf v})\;=\;[d(u_{n},v_{n})]\;\in\;\sss \N. \]
By the transfer principle, we have, for any three nonstandard
vertices $\bf u$, $\bf v$, and $\bf w$,
\[ {\bf d}({\bf u},{\bf v})\;\leq\;{\bf d}({\bf u},{\bf w})\,+\,{\bf d}({\bf w},{\bf v}). \]

We define the "galaxies" of $\sss D$ by first defining the "vertex
galaxies".  Two nonstandard vertices ${\bf u}=[u_{n}]$ and ${\bf
v}=[v_{n}]$ are defined to be in the same {\em vertex galaxy}
$\dot{\Gamma}$ of $\sss D$ if ${\bf d}({\bf u},{\bf v})$ is no
greater than a standard hypernatural $\bf k$, that is, if there
exists a natural number $k\in \N$ such that $\{n:d(u_{n},v_{n})\leq
k\}\in{\cal F}$.  In this case, we say that $\bf u$ and $\bf v$ are
{\em limitedly distant}, and we write ${\bf d}({\bf u},{\bf
v})\leq{\bf k}$.

By the same proof as that of \cite[Lemma 3.1]{gal}, we have

{\em Lemma 10.1.  The vertex galaxies partition the set of all
nonstandard vertices in $\sss D$.}

We define a galaxy $\Gamma$ of $\sss D$ as a vertex galaxy
$\dot{\Gamma}$ along with the set $A(\dot{\Gamma})$ of all the
nonstandard arcs that are each incident to two nonstandard vertices
in $\dot{\Gamma}$.  That is, for each $\dot{\Gamma}$, we have
$\Gamma =\{A(\dot{\Gamma}),\dot{\Gamma}\}$.  It follows from Lemma
10.1 that the galaxies of $\sss D$ partition $\sss D$ in the sense
that each nonstandard arc is in one and only one galaxy, namely, the
galaxy $\Gamma$ corresponding to the nonstandard vertex galaxy
$\dot{\Gamma}$.  We will say that all the nonstandard vertices in
$\dot{\Gamma}$ and all the nonstandard arcs  in $A(\dot{\Gamma})$
are in $\Gamma$.

The {\em principal galaxy} $\Gamma_{0}$ of $\sss D$ is that unique
galaxy, each of whose nonstandard vertices is limitedly distant from
some standard vertex.  All the vertices in $V$, where $V$ is the
vertex set in the standard digraph  $D$, are (i.e., can be
identified with) standard vertices in $\Gamma_{0}$, but there may be
other nonstandard vertices in $\Gamma_{0}$ as well.

Let us note that a galaxy need not be a subdigraph of $\sss D$
according to the definition of the latter adopted in Section 8 and
also in contrast to the terminology used in \cite{gal}.  This is
because a nonstandard vertex in $\dot{\Gamma}$ may have more
nonstandard tips in it than those belonging only to the nonstandard
arcs in $A(\dot{\Gamma})$.

Five examples of galaxies in nonstandard digraphs can be obtained
from Examples 3.2 to 3.6 in \cite{gal} simply by viewing the
branches therein as arcs.  This is because the galaxies in $\sss D$
are defined by means of distances based upon semipaths in the
digraph $D$, or equivalently by distances in the underlying graph
$G$ of $D$.

For the same reason, the theorems in \cite[Sections 3 and 4]{gal}
can be restated for digraphs with merely a change in wording.  Let
us list those appropriately reworded theorems here.  Their proofs
remain the same as those in \cite{gal}.  Remember that in this
section $D$ has an infinity of arcs, an infinity of vertices, and is
weakly connected.

{\em Theorem 10.2.  Let $D=\{A,V\}$ be locally finite (i.e., each
vertex in $V$ contains only finitely many tips).  Then, $\sss D$ has
at least one nonstandard vertex not in its principal galaxy
$\Gamma_{0}$, and thus at least one galaxy different from
$\Gamma_{0}$.}

Let $\Gamma_{a}$ and $\Gamma_{b}$ be two galaxies of $\sss D$ that
are different from the principal galaxy $\Gamma_{0}$ of $\sss D$. We
shall say that $\Gamma_{a}$ {\em is closer to $\Gamma_{0}$ than is
$\Gamma_{b}$} and that $\Gamma_{b}$ {\em is further away from
$\Gamma_{0}$ than is } $\Gamma_{a}$ if there are a ${\bf v}
=[v_{n}]$ in $\Gamma _{a}$ and a ${\bf w}=[w_{n}]$ in $\Gamma_{b}$
such that, for some ${\bf u}=[u_{n}]$ in $\Gamma_{0}$ and for every
$m\in \N$, we have
\[ N_{0}(m)\;=\;\{n:\,d(w_{n},u_{n})-d(v_{n},u_{n})\,\geq\, m\}\;\in\;{\cal F}. \]
Any set of galaxies for which every two of them, say $\Gamma_{a}$
and $\Gamma_{b}$ satisfy this condition will be said to be {\em
totally ordered according to their closeness to} $\Gamma_{0}$.  The
axioms for a total ordering using weak connectedness are easily
shown.  Also, these definitions do not depend upon the sequences
$\lla u_{n}\rra$, $\lla v_{n}\rra$, and $\lla w_{n}\rra$ chosen for
$\bf u$, $\bf v$ and $\bf w$.

{\em Theorem 10.3.  If $\sss D$ has a nonstandard vertex ${\bf v}$
that is not in its principal galaxy $\Gamma _{0}$, then there exists
a two-way infinite sequence of galaxies that is totally ordered
according to those galaxies closeness to $\Gamma _{0}$ with $\bf v$
being in one of those galaxies.}

The axioms of a partial ordering of a set of galaxies are same as
those for a total ordering except that the axiom of completeness is
dropped.

{\em Theorem 10.4.  Under the hypothesis of Theorem 10.3, the set of
all the galaxies of $\sss D$ is partially ordered according to the
closeness of the galaxies to the principal galaxy $\Gamma_{0}$.}

A subdigraph  $D_{s}$ of $D$ with the property that there is a
natural number $k$ such that $d(u,v)\leq k$ for all pairs of
vertices $u,v$ in $D_{s}$ will be called a {\em finitely dispersed}
subdigraph of $D$.  The structures of the galaxies other than the
principal galaxy $\Gamma_{0}$ are independent of any finitely
dispersed subdigraph of $D$ because the vertices $u_{n}$ in any
representative $\lla u_{n}\rra $ of any nonstandard vertex $\bf u$
in a galaxy other than $\Gamma_{0}$ must lie outside any finitely
dispersed subdigraph of $D$ for almost all $n$.  (This is the same
property as that for galaxies in nonstandard graphs \cite[Section 3
and Example 3.6]{gal}.)

\end{document}